\documentclass[12pt,twoside]{article}
\usepackage{amsmath}
\usepackage{amssymb,amsfonts}
\usepackage{graphicx} 
\usepackage{subcaption}                

\newcommand{\bH}{\mathbf{H}}

\newcommand{\bR}{\mathbf{R}}

\newcommand{\bS}{\mathbf{S}}

\newcommand{\cE}{\mathcal{E}}

\def\EE{\mathbf{E}}
\def\CC{\mathcal{C}}

\newcommand{\SXR}{\bS^2\!\times\!\bR}
\newcommand{\HXR}{\bH^2\!\times\!\bR}

\newcommand{\NIL}{\mathbf{Nil}}

\begin{document}


\title{Isoptic curves of cycloids}

\author{G\'eza~Csima \\
Budapest University of Technology and Economics, \\
Institute of Mathematics,
Department of Geometry Budapest, \\
P.O. Box 91, H-1521 \\
csgeza@math.bme.hu}
\maketitle
\thanks{}



\begin{abstract}
The history of the isoptic curves goes back to the 19th century, but nowadays the topic is experiencing a renaissance, providing numerous new results and new applications.
First, we define the notion of isoptic curve and outline some of the well-known results for strictly convex, closed curves. 
Overviewing the types of centered trochoids, we will be able to give the parametric equation of the isoptic curves of hypocycloids and epicycloids. Furthermore, we will determine the corresponding class of curves.
Simultaneously, we show that a generalized support function can be given to these types of curves in order to apply and extend the results for strictly convex, closed curves.

\end{abstract}

\newtheorem{theorem}{Theorem}[section]    
\newtheorem{claim}[theorem]{Claim}    
\newtheorem{lemma}[theorem]{Lemma}			
\newtheorem{corollary}[theorem]{Corollary}	
\newtheorem{definition}[theorem]{Definition}
\newtheorem{remark}[theorem]{Remark}
\newtheorem{conjecture}[theorem]{Conjecture}



\section{Introduction}

In this manuscript we work in the Euclidean plane $\EE^2$. Let us introduce the following definition:

\begin{definition}[\cite{yates}]The locus of the intersection of tangents
to a curve (or curves) meeting at a constant angle $\alpha$ $(0<\alpha <\pi)$ is
the $\alpha-$isoptic of the given curve (or curves). The isoptic curve with right angle called \textit{orthoptic curve}.
\label{defiso}
\end{definition}

Although the name "isoptic curve" was suggested by Taylor  in 1884 (\cite{T}), reference to former results can be found in \cite{yates}. In the obscure history of isoptic curves, we can find the names of la Hire (cycloids 1704) and Chasles (conics and epitrochoids 1837) among the contributors of the subject, however, the details of the research results are not available in English. A very interesting table of isoptic and orthoptic curves is introduced in \cite{yates}, unfortunately without any exact reference of its source. \emph{Our goal in this paper is to independently reconstruct some of the missing computations for the isoptic curves of hypocycloids and epicycloids and to extend the results presented in \cite{harom}, \cite{tizenketto} and \cite{picard}}.

However, recent works are available on the topic, which shows its timeliness. In \cite{harom} and \cite{tizenketto}, the Euclidean isoptic curves of closed strictly convex curves are studied using their support function.
Papers \cite{Kurusa,Wu71-1,Wu71-2} deal with Euclidean curves having a
circle or an ellipse for an isoptic curve. Further curves appearing as isoptic curves are well studied in Euclidean plane geometry $\mathbf{E}^2$, see e.g. \cite{Loria,Wi}.
Isoptic curves of conic sections have been studied in \cite{H} and \cite{Si}. There are results for Bezier curves by Kunkli et al. as well, see \cite{Kunkli}. Many papers focus on the properties of isoptics, e.g. \cite{nyolc,MM3,het}, and the references  therein. There are some generalizations of the isoptics as well \textit{e.g.} equioptic curves in \cite{Odehnal} by Odehnal or secantopics in \cite{tizenegy, tiz} by Skrzypiec. 

An algorithm for convex polyhedrons has been given by the authors in \cite{CsSz5} in order to generalize the notion of isoptic curve into the space and it has been developed by Kunkli et al. for non convex cases in \cite{Kunkli2}. The spatial case encompasses many applications in both physical and architectural aspects, see \cite{CsSz4}.

There are some results in non-Euclidean geometries as well. The isoptic curves of the hyperbolic line segment and proper conic sections are determined in \cite{CsSz1}, \cite{CsSz2} and \cite{CsSz3}. For generalized conic sections, and for their isoptics, see \cite{Cs-Sz20}. The isoptics of conic sections in elliptic geometry $\cE^2$ are determined in \cite{CsSz3}.

There are some results in three dimensional Thurson geometries as well. The isoptic surface of segments has been determined in \cite{CsSz-23} in $\NIL$ geometry and in \cite{Cs23} for $\SXR$ and $\HXR$ geometries. 

\section{Preliminary results}

In order to conduct further investigations on isoptics we need to summarize some preliminary results on the support function.

\begin{definition}
Let $\mathcal{C}$ be a closed, strictly convex curve which surrounds the origin. Let $p(t)$ where $t\in[0,2\pi[$ be the distance from $0$ to the support line of $\mathcal{C}$ being perpendicular to the vector $e^{it}$. The function $p$ is called a support function of $\mathcal{C}$.
\end{definition}

It is well-known \cite{egy} that the support function of a planar, closed, strictly convex curve $\mathcal{C}$ is differentiable. For now we would like to express the isoptic of $\CC$ using the support function. We claim the following lemma omitting the proof which can be found for example in \cite{ketto}.

\begin{lemma}[\cite{ketto}]
If $f(x,y,t)=0$ is a family of straight lines, then the equation of the envelope of these lines can be obtained by eliminating the variable $t$ from the two equations $f(x,y,t)=0$ and $\frac{d}{d t}f(x,y,t)=0$.
\label{l1}
\end{lemma}
This is used in \cite{ketto} to prove the following theorem.

\begin{theorem}[\cite{ketto}]
Given a planar, closed, strictly convex curve $\CC$ in polar coordinates with the radius $z$ a function of angle $t$, where $t\in\left[0,2\pi\right)$. Then the following equation holds \[z(t)=p(t)e^{it}+\dot{p}(t)i e^{it}.\]
\end{theorem}
The corollary of this theorem is that we may use this parametrization to determine the isoptic curve of $\CC$. The angle of $p(t)$ and $p(t+\pi-\alpha)$ is $\alpha$, since the $p(t)$, $p(t+\pi-\alpha)$ and their support lines determine a cyclic quadrilateral (see Figure 1.2). Our goal is to determine the intersection of these tangent lines which is the fourth vertex opposite the origin. A proof can be found in \cite{harom}.
\begin{theorem}[\cite{harom}]
Let $\CC$ be a plane, closed, strictly convex curve and suppose that the origin is in the interior of $\CC$. Let $p(t)$, $t\in\left[0,2\pi\right]$ be the support function of $\CC$. Then the $\alpha$-isoptic curve of $\CC$ has the form
\begin{equation}
z_{\alpha}(t)=p(t)e^{it}+\left(-p(t)\cot(\pi-\alpha)+\frac{1}{\sin(\pi-\alpha)}p(t+\pi-\alpha)\right)ie^{it}.
\label{eq:mozg}
\end{equation}
\label{theorem1}
\end{theorem}

\begin{definition}[\cite{yates}]
A \textbf{hypocycloid} is generated by a point on a circle rolling internally upon a fixed circle. An \textbf{epicycloid} is generated by a point on a circle rolling externally upon a fixed circle. A \textbf{hypotrochoid} is generated by a point rigidly attached to a circle rolling internally upon a fixed circle. An \textbf{epitrochoid} is generated by a point rigidly attached to a circle rolling externally upon a fixed circle. 
\end{definition}

We will use the following parametric equations of the hypo-and epicycloids, where we assumed that the radius of the fixed circle is $1$, and the radius of the rolling circle is rational $\frac{1}{a}:=\frac{p}{q}<1$ in its lowest terms, otherwise the curve never closes, and fills the space between the circles. Then we have exactly $p$ cusps and it is closed if and only if the length of parametric domain of $t$ is greater than equal to $2 q\pi$. In the case of hypocycloid, we also assume, that $2p\neq q$, which results in a segment.
\begin{equation}
\mathrm{Hypocycloid: }\left\{\frac{(a-1)\cos(t)+\cos((a-1) t)}{a},\frac{(a-1)\sin(t)-\sin ((a-1) t)}{a}\right\}
\label{eq:1}
\end{equation}
\begin{equation}
\mathrm{Epicycloid: }\left\{\frac{(a+1)\cos(t)-\cos((a+1) t)}{a},\frac{(a+1)\sin(t)-\sin ((a+1) t)}{a}\right\}
\label{eq:2}
\end{equation}

Finally, we need the parametric equations of the hypo-and epitrochoids:\\
Hypotrochoid:
\begin{equation}
\left\{(A-B)\cos(t)+H\cos\left(\frac{A-B}{B} t\right),(A-B)\sin(t)-H\sin \left(\frac{A-B}{B} t\right)\right\}
\label{eq:3}
\end{equation}
Epitrochoid:
\begin{equation}
\left\{(A+B)\cos(t)-H\cos\left(\frac{A+B}{B} t\right),(A+B)\sin(t)-H\sin \left(\frac{A+B}{B} t\right)\right\}
\label{eq:4}
\end{equation}
where the radius of the fixed and the rolling circles are $A$ and $B$ respectively, and $H$ is the distance of the rigid point to the center of the rolling circle (see \cite{Lawrence}).

\section{Isoptic curves}
Since the calculations of the isoptic curves of hypo-and epitrochoids are very similar, therefore we consider them together. Our first step, to determine the isoptic curves is always the tangent calculation. We need the derivative of the parametrization:
\begin{multline}
v_H(t)=\left\{-\frac{2 (a-1) \sin \left(\frac{a t}{2}\right) \cos \left(\frac{(a-2) t}{2} \right)}{a},
		\frac{2 (a-1) \sin \left(\frac{a t}{2}\right) \sin \left(\frac{(a-2) t}{2} \right) }{a}\right\}=\\
		=\frac{2 (a-1) \sin \left(\frac{a t}{2}\right)}{a}\left\{-\cos \left(\frac{(a-2) t}{2} \right), \sin \left(\frac{(a-2) t}{2} \right)  \right\}
\label{eq:5h}
\end{multline}
\begin{multline}
v_E(t)=\left\{\frac{2 (a+1) \sin \left(\frac{a t}{2}\right) \cos \left(\frac{(a+2) t}{2} \right)}{a},
		\frac{2 (a+1) \sin \left(\frac{a t}{2}\right) \sin \left(\frac{(a+2) t}{2} \right) }{a}\right\}=\\
		=\frac{2 (a+1) \sin \left(\frac{a t}{2}\right)}{a}\left\{\cos \left(\frac{(a+2) t}{2} \right), \sin \left(\frac{(a+2) t}{2} \right)  \right\}
\label{eq:5e}
\end{multline}
where we applied trigonometric product-to-sum and sum-to-product identities. 
\begin{remark} The tangent vector can be a null vector for discrete parameter values if $\sin \left(\frac{a t}{2}\right)=0$, but its direction may be determined in limit so that continuity remains.
\end{remark}
Now, it is easy to see, that the angle of two tangents is equal to the the angle of the corresponding tangent vectors. Considering the $t+\phi$ and $t-\phi$ parametric values:
\begin{equation}
\mathrm{H: }\ \ \frac{\left\langle v_H(t-\phi),v_H(t+\phi \right\rangle}{\left\|v_H(t-\phi)\right\|\left\|v_H(t+\phi)\right\|}=\cos((a-2)\phi)
\label{eq:6h}
\end{equation}
\begin{equation}
\mathrm{E: }\ \ \frac{\left\langle v_E(t-\phi),v_E(t+\phi \right\rangle}{\left\|v_E(t-\phi)\right\|\left\|v_E(t+\phi)\right\|}=\cos((a+2)\phi)
\label{eq:6e}
\end{equation}
that is \emph{independent form the parameter value} of $t$. This uniformity gives us the possibility to determine the isoptic curve. Let $\phi:=\frac{\alpha}{a\mp2}$ be true, if we are interested in the $\alpha$-isoptic curve. Then the angle of the oriented tangents that are drawn to points corresponding to the parameter values $t-\phi$ and $t+\phi$ is $\alpha$. 
\begin{remark}
In the case of the astroid ($a=4$), the value of $\phi$ is $\frac{\alpha}{2}$ so that the difference of considered two points in the parameter domain is exactly $\alpha$. 
\end{remark} 
From formulas(\ref{eq:5h}) and (\ref{eq:5e}), we can derive the equation of the tangent respected to the parameter $t$:
\begin{equation}
\mathrm{H:}\ \ x \sin \left(\frac{(a-2) t}{2} \right)+y \cos \left(\frac{(a-2) t}{2} \right)=\frac{(a-2) \sin \left(\frac{a t}{2}\right)}{a}
\label{eq:7h}
\end{equation}
\begin{equation}
\mathrm{E: }\ \ x \sin \left(\frac{(a+2) t}{2} \right)-y \cos \left(\frac{(a+2) t}{2} \right)=\frac{(a+2) \sin \left(\frac{a t}{2}\right)}{a}
\label{eq:7e}
\end{equation}
By replacing $t$ with $t-\phi$ and $t+\phi$, we get an equation system. We are looking for the common point of the above tangents that will be a point of the isoptic curve related to the parameters  $t$ and $\alpha$. Omitting the solution process and the simplification, here is it's result, which will be the parametrization of the isoptic curve, as well:
\begin{equation}
\mathrm{H: }\ \ x(t)=\frac{(a-2)  \left( \sin \left(\frac{(a-1) \alpha }{a-2}\right)\right)\cos (t)+\sin \left(\frac{\alpha }{a-2}\right) \cos ((a-1) t)}{a \sin (\alpha)}
\label{eq:parh1}
\end{equation}
\begin{equation}
\mathrm{H: }\ \ y(t)=\frac{(a-2)  \left(\sin \left(\frac{(a-1) \alpha }{a-2}\right)\sin (t) -\sin \left(\frac{\alpha }{a-2}\right) \sin ((a-1) t)\right)}{a \sin (\alpha)}
\label{eq:parh2}
\end{equation}
\begin{equation}
\mathrm{E: }\ \ x(t)=\frac{(a+2)  \left( \sin \left(\frac{(a+1) \alpha }{a+2}\right)\right)\cos (t)-\sin \left(\frac{\alpha }{a+2}\right) \cos ((a+1) t)}{a \sin (\alpha)}
\label{eq:pare1}
\end{equation}
\begin{equation}
\mathrm{E: }\ \ y(t)=\frac{(a+2)  \left(\sin \left(\frac{(a+1) \alpha }{a+2}\right)\sin (t) -\sin \left(\frac{\alpha }{a+2}\right) \sin ((a+1) t)\right)}{a \sin (\alpha)}
\label{eq:pare2}
\end{equation}
We can propose the following theorem realizing the similarities to (\ref{eq:3}) and (\ref{eq:4}):
\begin{theorem} Let us be given a hypocycloid with its parametrization
	$$\left\{\frac{(a-1)\cos(t)+\cos((a-1) t)}{a},\frac{(a-1)\sin(t)-\sin ((a-1) t)}{a}\right\}$$ 
	where $a=\frac{q}{p}$ and $t\in [0,2q\pi]$ such that $p,q\in\mathbb{Z}^+\wedge p<q\wedge 2p\neq q$. Then the $\alpha$-isoptic curve of it is a hypotrochoid given by the parametrization
	$$\left\{(A-B)\cos(t)+H\cos\left(\frac{A-B}{B} t\right),(A-B)\sin(t)-H\sin \left(\frac{A-B}{B} t\right)\right\},$$ 
	where
	$$A=\frac{(a-2) \sin \left(\frac{(a-1) \alpha }{a-2}\right)}{(a-1)\sin(\alpha)}, \ B=\frac{(a-2) \sin \left(\frac{(a-1) \alpha }{a-2}\right)}{a (a-1)\sin(\alpha)}, \ H=\frac{(a-2)  \sin \left(\frac{\alpha }{a-2}\right)}{a \sin(\alpha)}.$$
	\label{th:1h}
	\end{theorem}
	\begin{figure}[htp]
	\centering
		\includegraphics[width=0.48\textwidth]{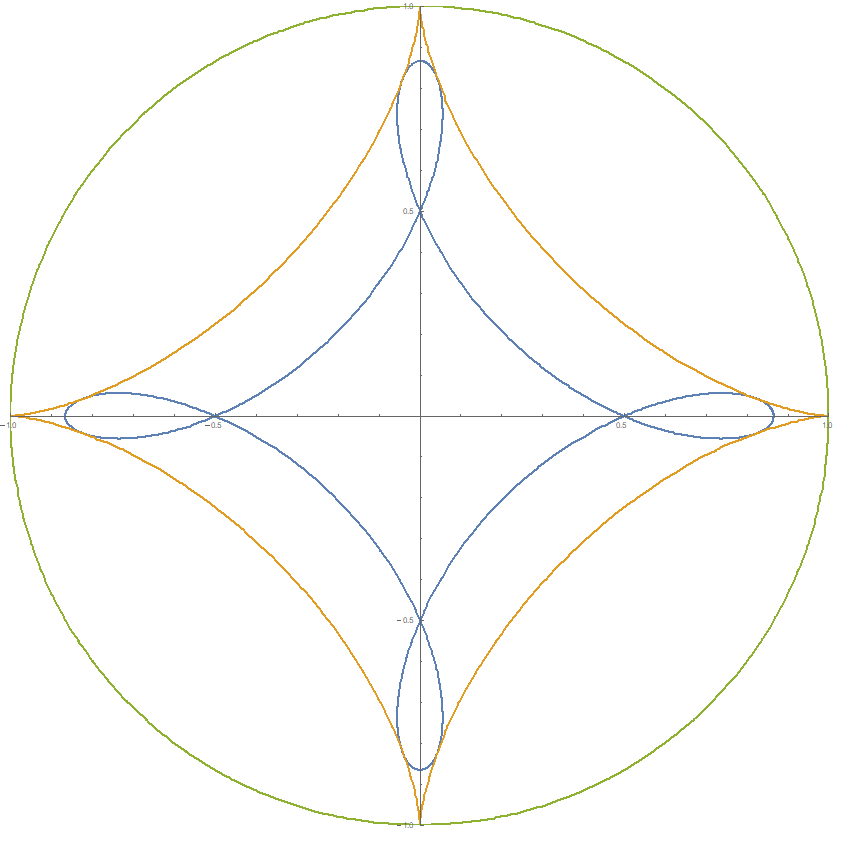} 		\includegraphics[width=0.48\textwidth]{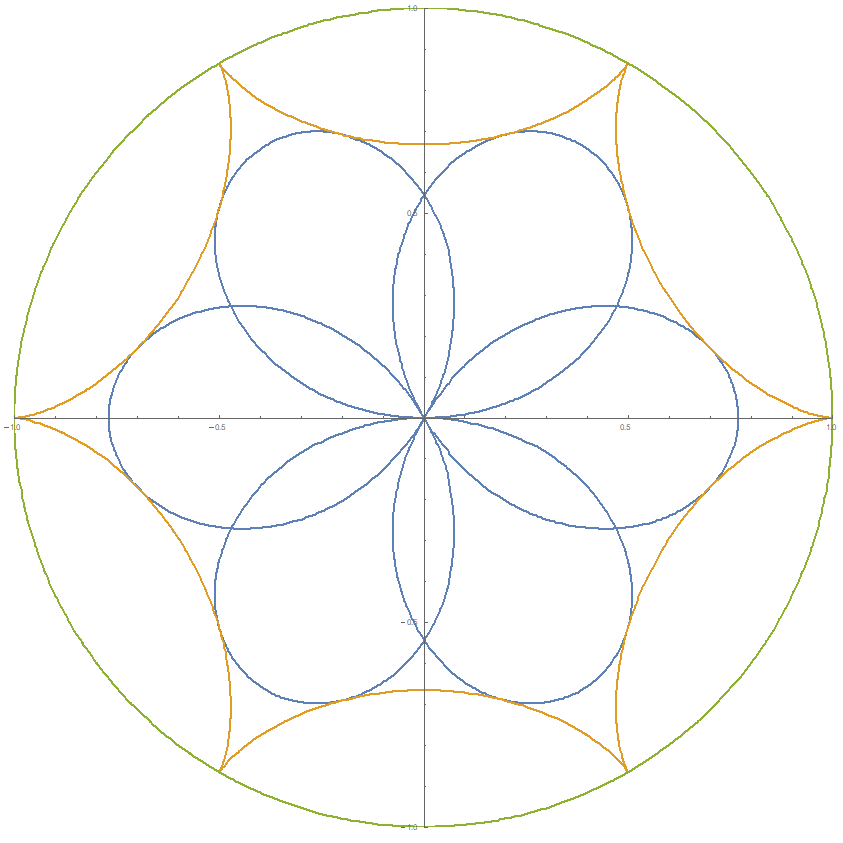}
	\caption{Isoptic curve for hypocycloid with $a=4,$ $\alpha=\frac{\pi}{3}$ (left) and $a=6,$ $\alpha=\frac{2\pi}{3}$ (right)}
	\label{fig:hypo}
\end{figure}

\begin{theorem} Let us be given an epicycloid with its parametrization
$$\left\{\frac{(a+1)\cos(t)-\cos((a+1) t)}{a},\frac{(a+1)\sin(t)-\sin ((a+1) t)}{a}\right\}$$ 
	where $a=\frac{q}{p}$ and $t\in [0,2q\pi]$ such that $p,q\in\mathbb{Z}^+\wedge p\leq q$. Then the $\alpha$-isoptic curve of it is an epitrochoid given by the parametrization
$$\left\{(A+B)\cos(t)-H\cos\left(\frac{A+B}{B} t\right),(A+B)\sin(t)-H\sin \left(\frac{A+B}{B} t\right)\right\},$$ 
	where
	$$A=\frac{(a+2) \sin \left(\frac{(a+1) \alpha }{a+2}\right)}{(a+1)\sin(\alpha)}, \ B=\frac{(a+2) \sin \left(\frac{(a+1) \alpha }{a+2}\right)}{a (a+1)\sin(\alpha)}, \ H=\frac{(a+2)  \sin \left(\frac{\alpha }{a+2}\right)}{a \sin(\alpha)}.$$
	\label{th:1e}
	\end{theorem}
	\begin{figure}[htp]
	\centering
		\includegraphics[width=0.48\textwidth]{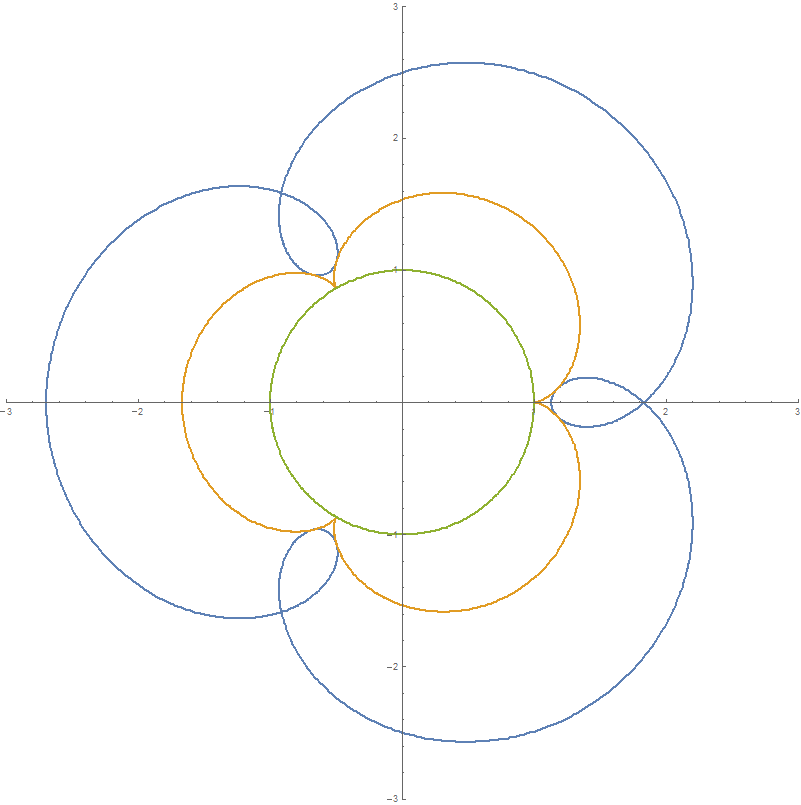} 		\includegraphics[width=0.48\textwidth]{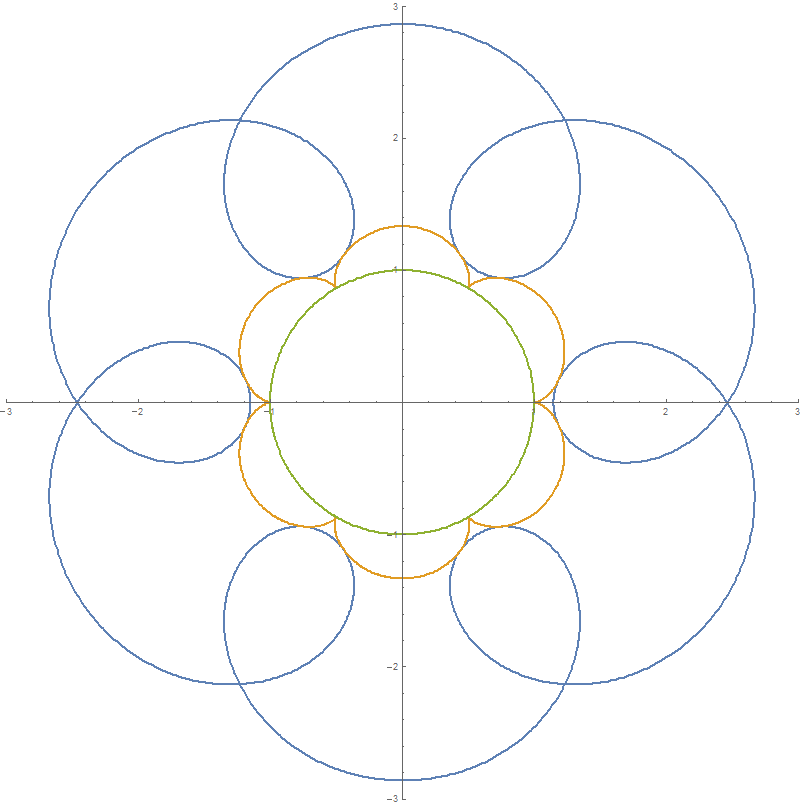}
	\caption{Isoptic curve for epicycloid with $a=3,$ $\alpha=\frac{\pi}{3}$ (left) and $a=6,$ $\alpha=\frac{\pi}{6}$ (right)}
	\label{fig:epi}
\end{figure}

\begin{remark}
It is easy to see, that for $\alpha=\dfrac{a-2}{a-1}\pi,$ $A=B=0$ in Theorem \ref{th:1h} and the resulted parametric curve is a circle, centered at the origin with radius  
$\dfrac{(a-2)  \sin \left(\frac{\pi }{a-1}\right)}{a \sin\left(\frac{a-2}{a-1}\right)}.$ For the epicycloid, in Theorem \ref{th:1e} $A=B=0$ if $\alpha=\dfrac{a+2}{a+1}\pi$ but that angle is greater than $\pi,$ therefore it is not a real isoptic curve.
\end{remark}
	\begin{figure}[htp]
	\centering
		\includegraphics[width=0.48\textwidth]{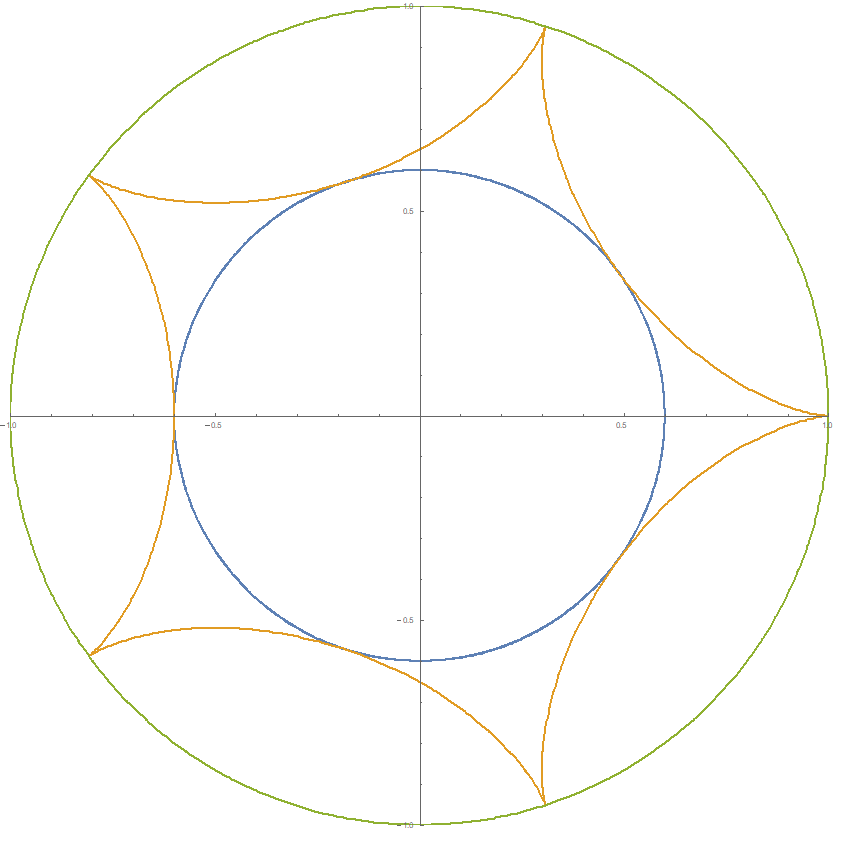} 
	\caption{Isoptic curve as a circle for hypocycloid with $a=5,$ $\alpha=\frac{3\pi}{4}$}
	\label{fig:hypo2}
\end{figure}

\section{Isoptic curves by support functions}	
One can realize that the tangents in formulas (\ref{eq:7h}) and (\ref{eq:7e}) are in Hesse form, therefore it is easy to calculate the distance of the line to the origin. Despite hypocycloids and epicycloids are non-convex curves, we can define their support function nonetheless in order to give another approach of the isoptic curve by Theorem \ref{theorem1}. We only have to apply a substitution: $t=\frac{2}{a-2}\left(\frac{\pi}{2}-u\right)$ in (\ref{eq:5h}) and $t=\frac{2}{a+2}\left(\frac{\pi}{2}-u\right)$ in (\ref{eq:5e}) to obtain:
\begin{equation}
\mathrm{H: }\ \ x \cos \left(u \right)+y \sin \left(u \right)=\frac{(a-2)}{a} \sin \left(\frac{a}{a-2}\left(\frac{\pi}{2}-u\right)\right),
\label{eq:8h}
\end{equation}
\begin{equation}
\mathrm{E: }\ \ x \cos \left(u \right)+y \sin \left(u \right)=\frac{(a+2)}{a} \sin \left(\frac{a}{a+2}\left(\frac{\pi}{2}-u\right)\right).
\label{eq:8e}
\end{equation}
It is easy to see, that the transverse vector of the tangent is $e^{iu}=\{\cos(u),\sin(u)\}$ and its distance to the origin is $\frac{(a-2)}{a} \sin \left(\frac{a}{a-2}\left(\frac{\pi}{2}-u\right)\right)$ in the case of the hypocycloid and $\frac{(a+2)}{a} \sin \left(\frac{a}{a+2}\left(\frac{\pi}{2}-u\right)\right)$ in the case of the epicycloid. Then we can define the quasi-support functions:
\begin{equation}
p_H(u)=\frac{(a-2)}{a} \sin \left(\frac{a}{a-2}\left(\frac{\pi}{2}-u\right)\right),
\label{eq:9h}
\end{equation}
\begin{equation}
p_E(u)=\frac{(a+2)}{a} \sin \left(\frac{a}{a+2}\left(\frac{\pi}{2}-u\right)\right).
\label{eq:9e}
\end{equation}
Now, we apply (\ref{eq:1}) from Theorem \ref{theorem1} to (\ref{eq:9h}) and to (\ref{eq:9e}) respectively:
\footnotesize
\[x_H(u)=
\frac{(a-2) \left(\sin (u) \left(\cot (\alpha ) \cos \left(\frac{\pi -a u}{a-2}\right)-\csc (\alpha ) \cos \left(\frac{\pi -a
   (\alpha +u)}{a-2}\right)\right)+\cos (u) \cos \left(\frac{\pi -a u}{a-2}\right)\right)}{a},\]
\[y_H(u)=\frac{(a-2) \left(\cos (u) \left(\csc (\alpha
   ) \cos \left(\frac{\pi -a (\alpha +u)}{a-2}\right)-\cot (\alpha ) \cos \left(\frac{\pi -a u}{a-2}\right)\right)+\sin (u) \cos
   \left(\frac{\pi -a u}{a-2}\right)\right)}{a},
\]
\[x_E(u)=
\frac{(a+2) \left(\sin (u) \left(\cot (\alpha ) \cos
   \left(\frac{a u+\pi }{a+2}\right)-\csc (\alpha ) \cos \left(\frac{a
   (\alpha +u)+\pi }{a+2}\right)\right)+\cos (u) \cos \left(\frac{a
   u+\pi }{a+2}\right)\right)}{a},\]
\[y_E(u)=\frac{(a+2) \left(\cos (u) \left(\csc
   (\alpha ) \cos \left(\frac{a (\alpha +u)+\pi }{a+2}\right)-\cot
   (\alpha ) \cos \left(\frac{a u+\pi }{a+2}\right)\right)+\sin (u) \cos
   \left(\frac{a u+\pi }{a+2}\right)\right)}{a}
\]
\normalsize
Since we are interested in the parametrization as the function of $t$, we take the inverse of the substitutions $u=\frac{1}{2}(\pi-t(a-2)),$ and $u=\frac{1}{2}(\pi-t(a+2))$ to obtain:
\footnotesize
\begin{equation}
x_H(t)=\frac{(a-2)  \left(\sin(\alpha +t)-\sin \left(t-\frac{a \alpha }{a-2}\right)+\sin \left(\frac{a \alpha }{a-2}-a t+t\right)-\sin (\alpha-a t +t)\right)}{2 a \sin (\alpha )}
\label{eq:parh3}
\end{equation}
\begin{equation}	
y_H(t)=\frac{(a-2) \left(-\cos (t-\alpha)+\cos \left(t+\frac{a
   \alpha }{a+2}\right)-\cos \left(a t+t-\frac{a \alpha }{a-2}\right)+\cos (a t +t-\alpha)\right)}{2 a \sin (\alpha )}
\label{eq:parh4}
\end{equation}

\begin{equation}
x_E(t)=\frac{(a+2)  \left(\sin(\alpha-t)+\sin \left(t+\frac{a \alpha }{a+2}\right)+\sin \left(\frac{a \alpha }{a+2}-a t-t\right)-\sin (\alpha-a t -t)\right)}{2 a \sin (\alpha )}
\label{eq:pare3}
\end{equation}
\begin{equation}	
y_E(t)=\frac{(a-2) \left(-\cos (\alpha +t)+\cos \left(t-\frac{a
   \alpha }{a-2}\right)-\cos \left(\frac{a \alpha }{a-2}-a t+t\right)+\cos (\alpha-a t +t)\right)}{2 a \sin (\alpha )}
\label{eq:pare4}
\end{equation}
\normalsize
We will show that the parametrization above results in the same curve as it has been described in (\ref{eq:parh1})--(\ref{eq:pare2}). 
Applying trigonometric sum-to-product identities for the first two, and second two parts of the numerators we obtain:
\begin{equation}
\mathrm{H:}\ \ 
x(t)=\frac{(a-2)  \left( \sin \left(\frac{(a-1) \alpha }{a-2}\right)\right)\cos (t-\frac{\alpha}{a-2})+\sin \left(\frac{\alpha }{a-2}\right) \cos ((a-1) (t-\frac{\alpha}{a-2}))}{a \sin (\alpha)}
\label{eq:parh5}
\end{equation}
\begin{equation}
\mathrm{H:}\ \ y(t)=\frac{(a-2)  \left(\sin \left(\frac{(a-1) \alpha }{a-2}\right)\sin (t-\frac{\alpha}{a-2}) -\sin \left(\frac{\alpha }{a-2}\right) \sin ((a-1) (t-\frac{\alpha}{a-2}))\right)}{a \sin (\alpha)}
\label{eq:parh6}
\end{equation}

\begin{equation}
\mathrm{E:}\ \ x(t)=\frac{(a+2)  \left( \sin \left(\frac{(a+1) \alpha }{a+2}\right)\right)\cos (\frac{\alpha}{a+2}-t)+\sin \left(\frac{\alpha }{a+2}\right) \cos ((a+1) (\frac{\alpha}{a+2}-t))}{a \sin (\alpha)}
\label{eq:pare5}
\end{equation}
\begin{equation}
\mathrm{E:}\ \ y(t)=\frac{(a+2)  \left(\sin \left(\frac{(a+1) \alpha }{a+2}\right)\sin (\frac{\alpha}{a+2}-t) -\sin \left(\frac{\alpha }{a+2}\right) \sin ((a+1) (\frac{\alpha}{a+2}-t))\right)}{a \sin (\alpha)}
\label{eq:pare6}
\end{equation}
Comparing (\ref{eq:parh1})--(\ref{eq:pare2}) to (\ref{eq:parh5})--(\ref{eq:pare6}), it can be easily verified that the two parametrization can be carried together. In the case of hypocycloid, we can arrange it by shifting the parameter domain with $\frac{\alpha}{a-2}$ but in the case of the epicycloid, we have to change the direction as well.

\begin{theorem}
Let us be given a $\CC$ hypocycloid with its parametrization
	$$\CC:\  \left\{\frac{(a-1)\cos(t)+\cos((a-1) t)}{a},\frac{(a-1)\sin(t)-\sin ((a-1) t)}{a}\right\}$$ 
	where $a=\frac{q}{p}$ and $t\in [0,2q\pi]$ such that $p,q\in\mathbb{Z}^+\wedge p<q\wedge 2p\neq q$.
 Let $p(t)$, $t\in\left[0,2\pi\right]$ be the support function of $\CC$. Then the $\alpha$-isoptic curve of $\CC$ has the form
\begin{equation}
z_{\alpha}(t)=p(t)e^{it}+\left(-p(t)\cot(\pi-\alpha)+\frac{1}{\sin(\pi-\alpha)}p(t+\pi-\alpha)\right)ie^{it},
\label{eq:mozg2h}
\end{equation}
where $p(t)=\frac{(a-2)}{a} \sin \left(\frac{a}{a-2}\left(\frac{\pi}{2}-t\right)\right).$
\label{theorem2e}
\end{theorem}

\begin{theorem}
Let us be given a $\CC$ epicycloid with its parametrization
	$$\CC:\ \ \left\{\frac{(a+1)\cos(t)-\cos((a+1) t)}{a},\frac{(a+1)\sin(t)-\sin ((a+1) t)}{a}\right\}$$ 
	where $a=\frac{q}{p}$ and $t\in [0,2q\pi]$ such that $p,q\in\mathbb{Z}^+\wedge p<q$.
 Let $p(t)$, $t\in\left[0,2\pi\right]$ be the support function of $\CC$. Then the $\alpha$-isoptic curve of $\CC$ has the form
\begin{equation}
z_{\alpha}(t)=p(t)e^{it}+\left(-p(t)\cot(\pi-\alpha)+\frac{1}{\sin(\pi-\alpha)}p(t+\pi-\alpha)\right)ie^{it},
\label{eq:mozg2e}
\end{equation}
where $p(t)=\frac{(a+2)}{a} \sin \left(\frac{a}{a+2}\left(\frac{\pi}{2}-t\right)\right).$
\label{theorem2h}
\end{theorem}


\end{document}